\documentclass[12pt]{amsart}
\usepackage{amsmath,graphicx,amssymb,amsthm}

\newtheorem{set}{set}[section]
\newtheorem{theorem}[set]{Theorem}
\theoremstyle{plain}

\newtheorem{corollary}[set]{Corollary}

\newtheorem{definition}[set]{Definition}

\newtheorem{lemma}[set]{Lemma}

\newtheorem{remark}[set]{Remark}

\begin{document}
\title[Haagerup's Approximation Property and Relative Amenability]
{Haagerup's Approximation Property and Relative Amenability}
\author{Jon P. Bannon}
\author{Junsheng Fang}
\address{Department of Mathematics, Siena College, Loudonville, NY 12211}
\email{jbannon@siena.edu}
\address{Department of Mathematics, University of New
Hampshire, Durham, NH 03824} \email{jfang@cisunix.unh.edu}

\begin{abstract}
A finite von Neumann algebra $\mathcal{M}$ with a faithful normal trace $%
\tau $ has Haagerup's approximation property (relative to a von
Neumann subalgebra $\mathcal{N}$) if there exists a net $(\varphi
_{\alpha })_{\alpha\in \Lambda}$ of normal completely positive
($\mathcal{N}$-bimodular) maps
from $\mathcal{M}$ to $\mathcal{M}$ that satisfy the subtracial condition $%
\tau \circ \varphi _{\alpha }\leq \tau $, the extension operators $%
T_{\varphi _{\alpha }}$ are bounded compact operators (in $\langle \mathcal{M%
},e_{\mathcal{N}}\rangle $), and pointwise approximate the identity in the
trace-norm, i.e., $\lim_{\alpha }||\varphi _{\alpha }(x)-x||_{2}=0$ for all $%
x\in \mathcal{M}$. We prove that the subtraciality condition can be
removed, and provide a description of Haagerup's approximation
property in terms of Connes's theory of correspondences. We show
that if $\mathcal{N}\subseteq
\mathcal{M}$ is an amenable inclusion of finite von Neumann algebras and $%
\mathcal{N}$ has Haagerup's approximation property, then
$\mathcal{M}$ also has Haagerup's approximation property. This work
answers two questions of Sorin Popa.
\end{abstract}

\maketitle

\section{Introduction}
A locally compact group $G$ has the \emph{Haagerup property} if
 there is a sequence of continuous
normalized  positive definite functions vanishing at infinity  on
$G$ that converges to $1$ uniformly on compact subsets of $G$.
In~\cite{Haa1}, Haagerup established the seminal result that free
groups have the Haagerup property. Now we know that the class of
groups having the Haagerup property is quite large. It includes all
compact groups, all locally compact amenable groups and all locally
compact groups that act properly on trees. There are many equivalent
characterizations of the Haagerup property. For instance, $G$ has
the Haagerup property if and only if there exists a continuous
positive real valued function $\psi$ on $G$ that is conditionally
negative definite and proper, i.e., $\lim_{g\rightarrow
\infty}\psi(g)=0$ (see~\cite{AW}). Also the Haagerup property is
equivalent to $G$ being a-T-menable in the sense of Gromov
(see~\cite{Gro1, Gro2, BCV}. An extensive treatment of the Haagerup
property for groups can be found in the book \cite{CCJJV}.
 Studying the class of
Haagerup groups has been a fertile endeavor. For example, the
Baum-Connes conjecture is solved for this class (see \cite{HigKa, Tu}).\\

The Haagerup property is a strong negation to Kazhdan's property T,
in that each of the equivalent definitions above stands opposite to
a definition of property T (see \cite{CCJJV}). A. Connes and V.
Jones introduced a notion of property T for von Neumann algebras in
terms of \emph{correspondences} \cite{CoJo}. Correspondences, as
introduced by Connes (see \cite{Co2, Co3, Po1}), are analogous to
group representations in the theory of von Neumann algebras. Connes
and Jones proved that a group has Kazhdan's property T if and only
if the associated group von Neumann algebra has the Connes-Jones
property T. Since property T has a natural analogue in the theory of
von Neumann algebras, we may expect that the same should be possible
for the Haagerup property. Such an analogue indeed exists, as Choda
proved in \cite{Cho} that a discrete group has the Haagerup property
if and only if its associated group von Neumann algebra has
Haagerup's approximation property first introduced in \cite{Haa1}:
There exists a net $(\varphi _{\alpha })_{\alpha
\in \Lambda }$ of normal completely positive maps from $\mathcal{M}$ to $%
\mathcal{M}$ such that (1) $\tau \circ
\varphi _{\alpha }(x^{\ast }x)\leq \tau (x^{\ast }x)$ for all $x\in \mathcal{%
M},$ (2)  the extension operator $T_{\varphi_\alpha}$ of
$\varphi_\alpha$ (see~\cite{Po2} or section 2.1 of this paper)
 is a compact operator in $\mathcal{B}(L^2(\mathcal{M},\tau))$ and (3) $\lim_{\alpha }||\varphi _{\alpha
}(x)-x||_{2}=0$ for all $x\in \mathcal{M}$. In~\cite{Po2}, Popa
asked if condition (1) can be removed for all finite von Neumann
algebras, and proved that if $\mathcal{M}$ is a non-$\Gamma$ type
${\rm II}_1$ factor, then (1) can be removed. In this paper, we
prove that (1) can be removed in general. \\

This enables us to provide a description of Haagerup's approximation
property in the language of correspondences. Suppose $\mathcal{M}$
is a finite von Neumann algebra with a faithful normal trace $\tau$
and $\mathcal{H}$ is a correspondence of $\mathcal{M}$. By the
Stinespring construction, if $\varphi$ is a normal completely
positive map from $\mathcal{M}$ to $\mathcal{M}$, there is a cyclic
correspondence $\mathcal{H}_\varphi$ of $\mathcal{M}$ associated to
$\varphi$. Every correspondence of $\mathcal{M}$ is equivalent to a
direct sum of cyclic correspondences associated to completely
positive maps as above (see~\cite{Po1}). We say that $\mathcal{H}$
is a \emph{${\rm C_0}$-correspondence} if $\mathcal{H}$ is
equivalent to $\oplus_{\alpha\in
\Lambda}\mathcal{H}_{\varphi_{\alpha}}$, where each
$\mathcal{H}_{\varphi_{\alpha}}$ is the correspondence of
$\mathcal{M}$ associated to a completely positive map
$\varphi_\alpha$ from $\mathcal{M}$ to $\mathcal{M}$  such that
 the extension operator $T_{\varphi_\alpha}$ of $\varphi_\alpha$
 is a compact operator in $\mathcal{B}(L^2(\mathcal{M},\tau))$.
 By the uniqueness of standard representation up to spatial
 isomorphism(see~\cite{Haa2}), the above definition of
 ${\rm C_0}$-correspondence does not depend on the choice of $\tau$.
  In this paper, we prove that $\mathcal{M}$ has
Haagerup's approximation property if and only if the identity
correspondence of $\mathcal{M}$ is weakly contained in some
 ${\rm C_0}$-correspondence of $\mathcal{M}$. We also show that if
  $\mathcal{M}$ has Haagerup's approximation property, then the equivalent class of ${\rm
 C_0}$-correspondences of $\mathcal{M}$ is dense in ${\rm
 Corr}(\mathcal{M})$, the set of equivalent classes of correspondences of
 $\mathcal{M}$. \\

In recent breakthrough work, Popa has combined relative versions of
property T and Haagerup's approximation property to create
\textquotedblleft deformation malleability\textquotedblright\
techniques to solve a number of old open questions about type ${\rm
II}_{1}$ factors (see~\cite{Po2, IPP}).  In  \cite{Po2}, Popa
 asked the following question: If $\mathcal{
N\subseteq M}$ is an amenable inclusion of finite von Neumann algebras and $\mathcal{N%
}$ has Haagerup's approximation property, does $\mathcal{M}$ also
have Haagerup's approximation property?  Recall that an inclusion
$\mathcal{N\subseteq M}$ of finite von Neumann algebras is amenable
if there exists a conditional expectation from the basic
construction $\langle \mathcal{M}, e_{\mathcal{N}}\rangle$ onto
$\mathcal{M}$ (see~\cite{Po1}). This question is motivated by the
analogous known result for groups:\ If $G$ is a  subgroup of a
discrete group $G_{0}$, and $G$ is co-F\o lner in $G_{0}$ in the
sense of Eymard~\cite{Eym}, then whenever $G$ has the Haagerup
property so does $G_{0}$. A proof of this result can be found in
\cite{CCJJV}. The condition that $G$ is co-F\o lner in $G_{0}$ in
the sense of Eymard is equivalent to the amenability of the
functorial inclusion $L(G)\subseteq L(G_{0})$ of
group von Neumann algebras \cite{Po2}.\\

In~\cite{Joli3}, Jolissaint proved that if the basic construction
$\langle \mathcal{M}, e_{{\mathcal N}}\rangle$ is a finite von
Neumann algebra and $\mathcal{N}$ has Haagerup's approximation
property, then $\mathcal{M}$ also has Haagerup's approximation
property. An affirmative answer to Popa's question for group von
Neumann algebras is found in the work of Anantharaman-Delaroche
\cite{AD2}. Anantharaman-Delaroche proved that the \emph{compact
approximation property} is equivalent to the Haagerup approximation
property in the group von Neumann algebra case. Recall that a
separable finite von Neumann algebra $\mathcal{M}$ has the compact
approximation property~\cite{AD2} if there exists a net $(\phi
_{\alpha })_{\alpha \in
\Lambda }$ of normal completely positive maps from $\mathcal{M}$ to $%
\mathcal{M}$ such that for all $x\in \mathcal{M}$ we have
$\lim_{\alpha }\phi _{\alpha }(x)=x$ ultraweakly and for all $\xi
\in L^{2}(\mathcal{M})$ and $\alpha \in \Lambda $, the map $x\mapsto
\phi _{\alpha }(x)\xi $ is a compact operator from \emph{the normed
space} $\mathcal{M}$ to $L^{2}(\mathcal{M})$. Anantharaman-Delaroche
proved that if $\mathcal{N\subseteq M}$ is an amenable inclusion and
$\mathcal{N}$ has the compact approximation property, then
$\mathcal{M}$ must have the compact approximation property. Using
properties (2) and (3) in the above definition of  Haagerup's
approximation property Anantharaman-Delaroche proved that Haagerup's
approximation property implies the compact approximation property.
It follows that if $\mathcal{N\subseteq M}$ is an amenable inclusion
and $\mathcal{N}$ has Haagerup's approximation property, then
$\mathcal{M}$ also has the compact approximation property. Popa's
question is then answered by appealing to the above result of Choda
to establish that for group von Neumann algebras the compact
approximation property
implies  Haagerup approximation property.\\

In this paper, we also answer Popa's second question affirmatively
for all finite von Neumann algebras. Our description of Haagerup's
approximation property in the language of correspondences plays a
key role in the proof.  The layout of the rest paper is as follows:%
\begin{equation*}
\begin{tabular}{l}
2. Preliminaries \\
3. Removal of the subtracial condition \\
4. ${\rm C_0}$-correspondences \\
5. Relative amenability and Haagerup's approximation property \\
\end{tabular}%
\end{equation*}

The first author wishes to express his deepest thanks to Paul Jolissaint and
Mingchu Gao for carefully reading an early version of this paper.

\section{Preliminaries}
\subsection{Extension of completely positive maps to Hilbert space operators}
Let $\mathcal{M}$ be a finite von Neumann algebra with a faithful
normal
trace $\tau$, and $\Omega_{\mathcal{M}}$ be the standard trace vector in $L^{2}(\mathcal{M%
},\tau )$ corresponding to $1\in \mathcal{M}$. For $x,y\in
\mathcal{M}$, $\langle x\Omega,y\Omega\rangle_\tau$ is defined to be
$\tau(y^*x)$ and $\|x\|_{2,\tau}=\tau(x^*x)^{1/2}$. When no
confusion arises, we simply write $\Omega$ instead of
$\Omega_{\mathcal{M}}$, and $\|x\|_{2}$ instead of
$\|x\|_{2,\tau}$.\\

 Suppose $\varphi$ is a normal completely positive map from
$\mathcal{M}$ to $\mathcal{M}$. Recall that if there
is a $c>0$ such that $\|\varphi(x)\|_2\leq c\|x\|_2$ for all $x\in \mathcal{M}$%
, then there is a (unique) bounded operator $T_\varphi$ on $L^2(\mathcal{M}%
,\tau )$ such that
\begin{equation*}
T_{\varphi }(x\Omega )=\varphi (x)\Omega \text{ \ \ \ \ }\forall
x\in \mathcal{M}\text{.}
\end{equation*}
$T_\varphi$ is called \emph{the extension operator} of $\varphi$.  If $%
\tau\circ \varphi\leq c_0\tau$ for some $c_0>0$, then
$\|\varphi(x)\|_2\leq c_0\|\varphi(1)\|^{1/2}$ $\|x\|_2$ (see Lemma
1.2.1 of~\cite{Po2}) and so there
is a bounded operator $T_\varphi$ on $L^2(\mathcal{M},\tau )$ such that $%
T_{\varphi }(x\Omega )=\varphi (x)\Omega$ for all $x\in
\mathcal{M}$.

\subsection{The basic construction and its compact ideal space}
Let $\mathcal{M}$ be a finite von Neumann algebra with a faithful
normal trace $\tau$, and $\mathcal{N}$ a von Neumann subalgebra of
$\mathcal{M}$. The basic construction $\langle
\mathcal{M},e_{\mathcal{N}}\rangle$ is the von Neumann algebra on
$L^2(\mathcal{M},\tau)$ generated by $\mathcal{M}$ and the
orthogonal projection $e_{\mathcal{N}}$ from $L^2(\mathcal{M},\tau)$
onto $L^2(\mathcal{N},\tau)$. Then $\langle
\mathcal{M},e_{\mathcal{N}}\rangle$ is a semi-finite von Neumann
algebra with a faithful normal semi-finite tracial weight ${\rm Tr}$
such that
\[{\rm Tr}(xe_{\mathcal{N}}y)=\tau(xy),\quad \forall x, y\in
\mathcal{M}.\] Recall that $\langle
\mathcal{M},e_{\mathcal{N}}\rangle=J\mathcal{N}'J$, where $J$ is the
conjugate linear isometry defined by $J(x\Omega)=x^*\Omega$,
$\forall x\in \mathcal{M}$. The compact ideal space of $\langle
\mathcal{M},e_{\mathcal{N}}\rangle$, denoted by $\mathcal{J}(\langle
\mathcal{M},e_{\mathcal{N}}\rangle)$, is the norm-closed two-sided
ideal generated by finite projections of $\langle
\mathcal{M},e_{\mathcal{N}}\rangle$. Since $e_{\mathcal{N}}$ is a
finite projection in $\langle \mathcal{M},e_{\mathcal{N}}\rangle$,
it follows that $e_{\mathcal{N}}\in \mathcal{J}(\langle
\mathcal{M},e_{\mathcal{N}}\rangle)$. We refer the reader
to~\cite{Jo,Po2} for more details on the basic construction and its
compact ideal space.

\subsection{Correspondences}
Let $\mathcal{N}$ and $\mathcal{M}$ be von Neumann algebras. A
 correspondence between $\mathcal{N}$ and $\mathcal{M}$ is a Hilbert space
$\mathcal{H}$ with a pair of commuting normal representations
$\pi_{\mathcal{N}}$ and $\pi_{\mathcal{M}^\circ}$ of  $\mathcal{N}$
and $\mathcal{M}^\circ$ (the opposite algebra of $\mathcal{M}$) on
$\mathcal{H}$, respectively. Usually, the triple $(\mathcal{H},
\pi_{\mathcal{N}}, \pi_{\mathcal{M}^\circ})$ will be denoted by
$\mathcal{H}$. For $x\in \mathcal{N}$, $y\in\mathcal{M}$ and
$\xi\in\mathcal{H}$, we shall write $x\xi y$ instead of
$\pi_{\mathcal{N}}(x)\pi_{\mathcal{M}^\circ}(y)\xi$. For two vectors
$\xi,\eta\in \mathcal{H}$, we denote by $\langle \xi,
\eta\rangle_{\mathcal{H}}$ the inner product of vectors $\xi$ and
$\eta$. If $\mathcal{N}=\mathcal{M}$, then we simply say
$\mathcal{H}$ is a correspondence of $\mathcal{M}$.\\

Two correspondences $\mathcal{H},\mathcal{K}$ between $\mathcal{N}$
and $\mathcal{M}$ are \emph{equivalent}, denoted by
$\mathcal{H}\cong\mathcal{K}$, if they are unitarily equivalent as
$\mathcal{N}-\mathcal{M}$ bimodules (see~\cite{Po1}).

\subsection{Correspondences associated to completely positive maps}
Let $\mathcal{M}$ be a finite von Neumann algebra with a faithful
normal trace $\tau$, and $\varphi$ be a normal completely positive
map from $\mathcal{M}$ to $\mathcal{M}$. Define on the linear space
$\mathcal{H}_0=\mathcal{M}\otimes \mathcal{M}$ the sesquilinear form
\[\langle x_1\otimes y_1, x_2\otimes
y_2\rangle_{\varphi}=\tau(\varphi(x_2^*x_1)y_1y_2^*),\quad \forall
x_1, x_2, y_1, y_2\in \mathcal{M}.
\]
It is easy to check that the complete positivity of $\varphi$ is
equivalent to the positivity of $\langle \cdot, \cdot
\rangle_{\varphi}$. Let $\mathcal{H}_\varphi$ be
 the completion  of $\mathcal{H}_0/\sim$,
  where $\sim$ is the equivalence  modulo the null space of
  $\langle\cdot, \cdot\rangle_{\varphi}$ in $\mathcal{H}_0$. Then $\mathcal{H}_\varphi$
  is a  correspondence of $\mathcal{M}$  and
the bimodule structure is given by $x(x_1\otimes y_1)y=xx_1\otimes
y_1 y$ (see~\cite{Po1}). We call $\mathcal{H}_\varphi$  the \emph
{correspondence of
$\mathcal{M}$ associated to} $\varphi$.\\

 The
correspondence $\mathcal{H}_{\rm id}$ associated to the identity
operator on $\mathcal{M}$ is called the \emph{identity
correspondence} of $\mathcal{M}$.  It is easy to see that
$\mathcal{H}_{\rm id}$ and $L^2(\mathcal{M},\tau)$ are equivalent as
correspondences of $\mathcal{M}$. The correspondence
$\mathcal{H}_{\rm co}$ associated to the rank one normal completely
positive map $\varphi(x)=\tau(x)1$ is called the \emph{coarse
correspondence} of $\mathcal{M}$.
 If $\mathcal{N}$ is a von
Neumann subalgebra of $\mathcal{M}$ and $E_{\mathcal{N}}$ is the
unique $\tau$-preserving normal conditional expectation from
$\mathcal{M}$ to $\mathcal{N}$, then the
 correspondence of $\mathcal{M}$ associated to $E_{\mathcal{N}}$ is
denoted by $\mathcal{H}_{\mathcal{N}}$ instead of
$\mathcal{H}_{E_{\mathcal{N}}}$.
\subsection{Left $\tau$-bounded vectors}
Let $\mathcal{N}, \mathcal{M}$ be  finite von Neumann algebras with
 faithful normal traces $\tau_\mathcal{N}$ and $\tau_\mathcal{M}$, respectively,
  and $\mathcal{H}$ be a correspondence between $\mathcal{N}$ and $\mathcal{M}$.
   Let $\xi\in \mathcal{H}$ be a vector. Recall that $\xi$ is a
\emph{left $($or right$)$ $\tau$-bounded vector} if there is a
positive number $K$ such that $\langle \xi,\xi
x\rangle_{\mathcal{H}}\leq K\tau_\mathcal{M}(x)$ (or $\langle
x\xi,\xi\rangle_{\mathcal{H}}\leq K\tau_\mathcal{N}(x)$,
respectively) for all $x\in \mathcal{N}_+$ (or $x\in \mathcal{M}_+$,
respectively). A vector $\xi$ is called a \emph{$\tau$-bounded
vector} if it is both left $\tau$-bounded and right $\tau$-bounded.
The set of $\tau$-bounded vectors is a dense vector subspace of
$\mathcal{H}$ (see Lemma 1.2.2 of~\cite{Po1}).
\subsection{Coefficients}
Let $\mathcal{N}, \mathcal{M}$ be  finite von Neumann algebras with
 faithful normal traces $\tau_\mathcal{N}$ and $\tau_\mathcal{M}$, respectively,
  and $\mathcal{H}$ be a correspondence between $\mathcal{N}$ and $\mathcal{M}$.
   For a left $\tau$-bounded
vector $\xi$, we can define a bounded operator $T:\,
L^2(\mathcal{M},\tau_\mathcal{M})\rightarrow \mathcal{H}$ by
$T(y\Omega_{\mathcal{M}})=\xi y$ for every $y\in\mathcal{M}$. Let
$\Phi_\xi(x)=T^*\pi_\mathcal{N}(x)T$, where $\pi_\mathcal{N}(x)$ is
the left action of $x\in \mathcal{N}$ on $\mathcal{H}$. Then
$\Phi_\xi$ is a normal completely positive map from $\mathcal{N}$ to
$\mathcal{M}$ (see 1.2.1 of~\cite{Po1}). $\Phi_\xi$ is called the
\emph{coefficient} corresponding to $\xi$, which is uniquely
determined by
\begin{equation}\label{Eq:Coefficient}
\langle
\Phi_\xi(x)y\Omega_{\mathcal{M}},z\Omega_{\mathcal{M}}\rangle_{\tau_\mathcal{M}}=\langle
x\xi y, \xi z\rangle_{\mathcal{H}}
\end{equation}
for all $x\in \mathcal{N}$ and $y,z\in \mathcal{M}$. Therefore, \[
\Phi_\xi(x)=\frac{d\langle
x\xi\cdot,\xi\rangle_{\mathcal{H}}}{d\tau_\mathcal{M}},\, {\rm
i.e}.,\, \tau_\mathcal{M}(\Phi_\xi(x)y)=\langle x\xi
y,\xi\rangle_{\mathcal{H}},\,\forall x\in \mathcal{N}, y\in
\mathcal{M}.\]
 If $\mathcal{N}=\mathcal{M}$,
$\tau_\mathcal{N}=\tau_\mathcal{M}$, and $x\geq 0$,
\[\tau_\mathcal{M}(\Phi_\xi(x))=\langle \Phi_\xi(x)\Omega_{\mathcal{M}},\Omega_{\mathcal{M}}\rangle_\tau=\langle x\xi, \xi
\rangle_{\mathcal{H}}\leq K\tau_\mathcal{M}(x).
\]
 By Lemma 1.2.1 of~\cite{Po2},
$\Phi_\xi$ can be extended to a bounded operator $T_{\Phi_\xi}$ from
$L^2(\mathcal{M},\tau)$ to $L^2(\mathcal{M},\tau)$.\\

It follows by a maximality argument that $\mathcal{H}$ is a direct
sum of cyclic correspondences associated to coefficients as above.
\subsection{Composition of correspondences} Suppose $\mathcal{M},
\mathcal{N}, \mathcal{P}$ are finite von Neumann algebras, and
$\tau_{\mathcal{P}}$ is a faithful normal trace on $\mathcal{P}$.
 Let
$\mathcal{H}$ be a  correspondence between $\mathcal{N}$ and
$\mathcal{P}$ and $\mathcal{K}$ be a  correspondence between
$\mathcal{P}$ and $\mathcal{M}$. Let $\mathcal{H}'$ and
$\mathcal{K}'$ be vector subspaces of the $\tau$-bounded vectors in
$\mathcal{H}$ and $\mathcal{K}$, respectively. For $\xi_1,\xi_2\in
\mathcal{H}'$ and $\eta_1,\eta_2\in\mathcal{K}'$,
\[\langle \xi_1\otimes \eta_1, \xi_2\otimes \eta_2\rangle=
\langle\xi_1 p, \xi_2 \rangle_{\mathcal{H}}=\langle
q\eta_1,\eta_2\rangle_{\mathcal{K}}=\tau_{\mathcal{P}}(qp)
\]
defines an inner product on $\mathcal{H}'\otimes \mathcal{K}'$,
where $p$ and $q$ are  Radon-Nikodym derivatives of  normal linear
forms $\mathcal{P}\ni z\rightarrow \langle
z\eta_1,\eta_2\rangle_{\mathcal{K}}$ and $\mathcal{P}\ni
z\rightarrow \langle\xi_1 z, \xi_2 \rangle_{\mathcal{H}}$ with
respect to the trace $\tau_{\mathcal{P}}$, respectively
(see~\cite{Po1}). The \emph{composition correspondence} (or the
\emph{tensor product correspondence})
$\mathcal{H}\underset{\mathcal{P}}{\otimes}\mathcal{K}$ is the
completion of $\mathcal{H}'\otimes \mathcal{K}'/\sim$, where $\sim$
is the equivalence  modulo the null space of
  $\langle\cdot, \cdot\rangle$ in $\mathcal{H}'\otimes \mathcal{K}'$, and
the $\mathcal{N}-\mathcal{M}$ bimodule structure is given by
$x(\xi\otimes \eta)y=x\xi\otimes \eta y$. It is easy to verify that
the composition of correspondences is associative.
\subsection{Induced correspondences}
A very important operation in various representation theories (e.g.
for groups) is that of inducing from smaller objects to larger ones.
We also have such a concept equally important to the theory of
correspondences. Let $\mathcal{M}$ be a finite von Neumann algebra
with a faithful normal trace $\tau$ and $\mathcal{N}$
  a  von Neumann subalgebra. If
$\mathcal{H}$ is a correspondence of $\mathcal{N}$, then the
\emph{induced correspondence} by $\mathcal{H}$ from $\mathcal{N}$ up
to $\mathcal{M}$ is ${\rm
Ind}_{\mathcal{N}}^{\mathcal{M}}(\mathcal{H})=L^2(\mathcal{M},\tau)\underset{\mathcal{N}}\otimes
\mathcal{H} \underset{\mathcal{N}}\otimes L^2(\mathcal{M},\tau)$,
where the first $L^2(\mathcal{M})$ is regarded as a left
$\mathcal{M}$ and right $\mathcal{N}$ module and the second
$L^2(\mathcal{M})$ is regarded as a left $\mathcal{N}$ and right
$\mathcal{M}$ module. If $\mathcal{H}$ is the identity
correspondence of $\mathcal{N}$, then ${\rm
Ind}_{\mathcal{N}}^{\mathcal{M}}(\mathcal{H})$ is the correspondence
$\mathcal{H}_{\mathcal{N}}$ of $\mathcal{M}$(see Proposition 1.3.6
of~\cite{Po1}).

\subsection{Relative amenability} Let $\mathcal{H},\mathcal{K}$ be
two  correspondences between $\mathcal{N}$ and $\mathcal{M}$. We say
that $\mathcal{H}$ is \emph{weakly contained} in $\mathcal{K}$,
 if for every
$\epsilon>0$, and finite subsets $E\subseteq \mathcal{N}$,
$F\subseteq\mathcal{M}$, $\{\xi_1,\cdots,\xi_n\}\subseteq
\mathcal{H}$,  there exists
$\{\eta_1,\cdots,\eta_n\}\subseteq\mathcal{K}$ such that
\[|\langle x\xi_i y,\xi_j\rangle_{\mathcal H}-\langle x\eta_i y,\eta_j\rangle_{\mathcal
H}|<\epsilon,
\] for all $x\in E, y\in F$ and $1\leq i, j\leq n$. If
$\mathcal{H}$ is weakly contained in $\mathcal{K}$, we will  denote
this  by $\mathcal{H}\prec \mathcal{K}$. We refer the reader
to~\cite{CoJo,Po1} for more details on weak containment and topology
on correspondences. \\

Let $\mathcal{M}$ be a finite von Neumann algebra with a faithful
normal trace $\tau$, and $\mathcal{N}$ a von Neumann subalgebra.
Recall that $\mathcal{M}$ is \emph{relative amenable to}
$\mathcal{N}$ if
$\mathcal{H}_{\rm{id}}\prec\mathcal{H}_{\mathcal{N}}$. The algebra
$\mathcal{M}$ is relative amenable to $\mathcal{N}$ if and only if
there exists a conditional expectation from the basic construction
$\langle \mathcal{M}, e_{\mathcal{N}}\rangle$ onto $\mathcal{M}$
(see~\cite{Po1}).
\section{Removal of the Subtracial Condition}

The following definition is given by Popa in~\cite{Po2}.

\begin{definition}
\label{HAP}\emph{\ Let $\mathcal{M}$ be a finite von Neumann algebra and $%
\mathcal{N}$ a von Neumann subalgebra. $\mathcal{M}$ has  \emph{%
Haagerup's approximation property relative to} $\mathcal{N}$ if
there exists
a normal faithful trace $\tau$ on $\mathcal{M}$ and a net of normal completely positive $%
\mathcal{N}$-bimodular maps $\{\varphi _{\alpha }\}_{\alpha \in \Lambda }$
from $\mathcal{M}$ to $\mathcal{M}$ satisfying the conditions:
\begin{enumerate}
\item[$1.$]\quad  $\tau \circ \varphi _{\alpha }\leq \tau $, $\forall\alpha \in \Lambda
$;
\item[$2.$]\quad
 $T_{\varphi _{\alpha }}\in \mathcal{J}(\langle\mathcal{M},e_\mathcal{N}\rangle)$,  $\forall \alpha \in \Lambda
$;
\item[$3.$]\quad  $\lim_{\alpha }\Vert \varphi _{\alpha }(x)-x\Vert _{2}=0$,  $
\forall x\in \mathcal{M}$.
\end{enumerate} }
\end{definition}

In~\cite{Po2}(Remark 2.6), Popa asked  if the condition $1$ in
Definition~\ref{HAP} can be removed or not (see Remark 2.6 of
\cite{Po1}). In this section we give an affirmative answer to Popa's
question. Precisely, we will prove the following theorem.

\begin{theorem}
\label{Fang1} Let $\mathcal{M}$ be a finite von Neumann algebra with a
faithful normal trace $\tau$ and $\mathcal{N}$ a von Neumann subalgebra.
Suppose $\{\varphi_\alpha\}_{\alpha\in \Lambda}$ is a net of normal completely positive $%
\mathcal{N}$-bimodular  maps from $\mathcal{M}$ to $\mathcal{M}$
satisfying the conditions $2$ and $3$ as in definition~\ref{HAP},
i.e. $\lim_{\alpha }\Vert \varphi _{\alpha }(x)-x\Vert _{2}=0$ for
all $x\in \mathcal{M}$ and the map $x\Omega\rightarrow
\varphi_\alpha(x)\Omega$ extends to a
bounded operator $T_{\varphi_\alpha}$ in $\mathcal{J}(\langle\mathcal{M},e_%
\mathcal{N}\rangle)$ for every $\alpha\in \Lambda$. Then there exists a net $%
\{\psi_\beta\}_{\beta\in \Gamma}$ of normal completely positive
$\mathcal{N}$-bimodular maps from $\mathcal{M}$ to $\mathcal{M}$
satisfying

\begin{enumerate}
\item[$1^{\prime }.$] \quad $\psi_\beta(1)=1$ and $\tau\circ \psi_\beta=\tau$%
, $\forall \beta\in \Gamma$;

\item[$2^{\prime }.$] \quad $T_{\psi_\beta}\in \mathcal{J}(\langle\mathcal{M}%
,e_\mathcal{N}\rangle)$,  $\forall \beta \in \Gamma$;

\item[$3^{\prime }.$] \quad $\lim_\beta\|\psi_\beta(x)-x\|_{2}=0$, $%
\forall x\in \mathcal{M}$.
\end{enumerate}

In particular, $\mathcal{M}$ has Haagerup's approximation property
relative to $\mathcal{N}$.
\end{theorem}

The ideas of the proof of Theorem~\ref{Fang1} are from Lemma 1.1.1 of \cite%
{Po2}, Day's trick \cite{Day}, and Proposition 2.1 of~\cite{Joli3}.
The following lemma is 2 of Lemma 1.1.1 of \cite{Po2} with a minor
change. For the sake of completeness, we include the proof.

\begin{lemma}
\label{Po1} Let $\varphi$ be a normal completely positive $\mathcal{N}$-bimodular map from $%
\mathcal{M}$ to $\mathcal{M}$. Let $a=1\vee \varphi(1)$ and $%
\varphi^{\prime}(\cdot)=a^{-1/2}\varphi(\cdot)a^{-1/2}$. Then
$\varphi^{\prime }$
is a normal completely positive $\mathcal{N}$-bimodular map from $\mathcal{M}$ to $\mathcal{M%
}$ and satisfies $\varphi^{\prime }(1)\leq 1$, $\tau\circ\varphi^{\prime
}\leq \tau\circ\varphi$ and the estimate:
\begin{equation*}
\|\varphi^{\prime }(x)-x\|_{2}\leq
\|\varphi(x)-x\|_2+2\|\varphi(1)-1\|_2^{1/2}\cdot\|x\|,\quad\forall x\in
\mathcal{M}.
\end{equation*}
\end{lemma}

\begin{proof}
Since $a\in \mathcal{N}^{\prime }\cap \mathcal{M}$, $\varphi^{\prime }$ is $%
\mathcal{N}$-bimodular.  We clearly have $\varphi ^{\prime }(1)\leq
1$. Since $a^{-1}\leq 1$, for $x\geq 0$ we get $\tau (\varphi
^{\prime}(x))\leq \tau (\varphi (x))$. Also, we have:
\begin{equation*}
\Vert \varphi ^{\prime }(x)-x\Vert _{2}\leq \Vert a^{-1/2}(\varphi
(x)-x)a^{-1/2}\Vert _{2}+\Vert a^{-1/2}xa^{-1/2}-x\Vert _{2}
\end{equation*}%
\begin{equation*}
\leq \Vert \varphi (x)-x\Vert _{2}+2\Vert a^{-1/2}-1\Vert _{2}\cdot \Vert
x\Vert .
\end{equation*}%
By the Powers-St\o rmer inequality (also see Proposition 1.2.1 of \cite{Co1}%
),
\begin{equation*}
\Vert a^{-1/2}-1\Vert _{2}\leq \Vert a^{-1}-1\Vert _{1}^{1/2}=\Vert
a^{-1}-aa^{-1}\Vert _{1}^{1/2}
\end{equation*}%
\begin{equation*}
\leq \Vert a-1\Vert _{2}^{1/2}\Vert a^{-1}\Vert _{2}^{1/2}\leq \Vert \varphi
(1)-1\Vert _{2}^{1/2}.
\end{equation*}%
Thus,
\begin{equation*}
\Vert \varphi ^{\prime }(x)-x\Vert _{2}\leq \Vert \varphi (x)-x\Vert
_{2}+2\Vert \varphi (1)-1\Vert _{2}^{1/2}\cdot \Vert x\Vert .
\end{equation*}
\end{proof}

\begin{lemma}
\label{Po2} Let $\varphi$ be a normal completely positive $\mathcal{N}$-bimodular map from $%
\mathcal{M}$ to $\mathcal{M}$ such that $\varphi(1)\leq 1$. Let $b=1\vee
(d\tau\circ\phi/d\tau)$ and $\varphi^{\prime}(\cdot)= \varphi(b^{-1/2}\cdot
b^{-1/2})$. Then $\varphi^{\prime }$ is a normal  completely positive $\mathcal{N}$%
-bimodular map from $\mathcal{M}$ to $\mathcal{M}$ and satisfies $%
\varphi^{\prime }(1)\leq \varphi(1)\leq 1$, $\tau\circ\varphi^{\prime }\leq
\tau$ and the estimate:
\begin{equation*}
\|\varphi^{\prime }(x)-x\|_{2}^2\leq
2\|\varphi(x)-x\|_2+5\|\tau\circ\phi-\tau\|^{1/2}\cdot\|x\|,\quad\forall
x\in \mathcal{M}.
\end{equation*}
\end{lemma}

\begin{proof}
Note that $\Vert \tau \circ \varphi -\tau \Vert=\Vert d\tau \circ \varphi
/d\tau -1\Vert _{1}$ and $\Vert b-1\Vert _{1}\leq \Vert d\tau \circ \varphi
/d\tau -1\Vert _{1}$. Now Lemma~\ref{Po2} follows simply from 3 of Lemma
1.1.1 of \cite{Po2}.
\end{proof}

\begin{lemma}
\label{Derivative}Let $\varphi$ be a normal completely positive
$\mathcal{N}$-bimodular map
from $\mathcal{M}$ to $\mathcal{M}$ such that $\varphi(1)\leq 1$ and $%
\tau\circ\varphi\leq \tau$. Let $h=\varphi(1)$ and $k=d\tau\circ\phi/d\tau$.
Then $0\leq h,k\leq 1$, $h,k\in \mathcal{N}^{\prime }\cap\mathcal{M}$, and $%
E_{\mathcal{N}}(h)=E_{\mathcal{N}}(k)$.
\end{lemma}

\begin{proof}
It is easy to see that $0\leq h,k\leq 1$ and $h,k\in \mathcal{M}$. Since $%
\varphi$ is $\mathcal{N}$-bimodular,
\begin{equation*}
bh=b\varphi(1)=\varphi(b)=\varphi(1)b=hb,\quad\forall b\in\mathcal{N}.
\end{equation*}
Note that for all $x\in \mathcal{M}$ and $b\in \mathcal{N}$,
\begin{equation*}
\tau(x(bk-kb))=\tau(xbk)-\tau(bxk)=\tau(\phi(xb))-\tau(\phi(bx))
\end{equation*}%
\begin{equation*}
=\tau(\phi(x)b)-\tau(b\phi(x))=0,
\end{equation*}
and
\begin{equation*}
\tau(E_{\mathcal{N}}(h)x)=\tau(hE_{\mathcal{N}}(x))=\tau(\varphi(1)E_{%
\mathcal{N}}(x))=\tau(\varphi(E_{\mathcal{N}}(x)))
\end{equation*}
\begin{equation*}
=\tau(E_{\mathcal{N}}(x)k)=\tau(xE_{\mathcal{N}}(k))=\tau(E_{\mathcal{N}%
}(k)x).
\end{equation*}
Hence, $bh=hb$ and $E_{\mathcal{N}}(h)=E_{\mathcal{N}}(k)$.
\end{proof}

\begin{lemma}
\label{equation}Let $\varphi$ be a normal completely positive
$\mathcal{N}$-bimodular map from $\mathcal{M}$ to $\mathcal{M}$ such
that $\varphi(1)\leq 1-\epsilon $ for some $\epsilon>0$ and
$\tau\circ\varphi\leq \tau$. Let $h=\varphi(1)$ and
$k=d\tau\circ\phi/d\tau$. Then there exist positive operators
$a,b\in \mathcal{N}^{\prime }\cap \mathcal{M}$ such that
\begin{equation*}
1-h=aE_{\mathcal{N}}(b)\quad\text{and}\quad 1-k=E_{\mathcal{N}}(a)b.
\end{equation*}
\end{lemma}

\begin{proof}
Let $b=1-k$. By Lemma~\ref{Derivative}, $b$ is a positive operator in $%
\mathcal{N}^{\prime }\cap \mathcal{M}$ and
\begin{equation*}
E_{\mathcal{N}}(b)=1-E_{\mathcal{N}}(k)=1-E_{\mathcal{N}}(h)=E_{\mathcal{N}%
}(1-h).
\end{equation*}
Since $h\leq 1-\epsilon$, $1-h\geq \epsilon$ and therefore $E_{\mathcal{N}%
}(1-h)\geq \epsilon$. Hence $(E_{\mathcal{N}}(1-h))^{-1}$ exists. For all $%
b\in\mathcal{N}$, by Lemma~\ref{Derivative},
\begin{equation*}
bE_{\mathcal{N}}(1-h)=E_{\mathcal{N}}(b(1-h))=E_{\mathcal{N}}((1-h)b)=E_{%
\mathcal{N}}(1-h)b.
\end{equation*}%
Hence, $E_{\mathcal{N}}(1-h)\in \mathcal{N}\cap \mathcal{N}^{\prime }$ and $%
(E_{\mathcal{N}}(1-h))^{-1}\in \mathcal{N}\cap \mathcal{N}^{\prime }$. So $%
a=(1-h)(E_{\mathcal{N}}(1-h))^{-1}$ is a positive operator in $\mathcal{N}%
^{\prime }\cap \mathcal{M}$. Since $E_{\mathcal{N}}(b)=E_{\mathcal{N}}(1-h)$%
, it is routine to check that $1-h=aE_{\mathcal{N}}(b)$ and $1-k=E_{\mathcal{%
N}}(a)b.$
\end{proof}

\begin{proof}[Proof of Theorem~\protect\ref{Fang1}]
Let $a_{\alpha }=1\vee \varphi _{\alpha }(1)$ and $\varphi _{\alpha
}^{\prime }(\cdot )=a_{\alpha }^{-1/2}\varphi _{\alpha }(\cdot )a_{\alpha
}^{-1/2}$. By Lemma~\ref{Po1}, $\{\varphi _{\alpha }^{\prime }\}_{\alpha }$
satisfy the condition $3^{\prime }$ in Theorem~\ref{Fang1} and $\varphi
_{\alpha }^{\prime }(1)\leq 1$ for every $\alpha \in \Lambda $. By Lemma~\ref{Derivative}, $%
T_{\varphi _{\alpha }^{\prime }}=a_{\alpha }^{-1/2}Ja_{\alpha
}^{-1/2}JT_{\varphi _{\alpha }}\in
\langle\mathcal{M},e_\mathcal{N}\rangle$. Since
  $T_{\varphi _{\alpha }}\in \mathcal{J}%
(\langle\mathcal{M},e_\mathcal{N}\rangle)$, $T_{\varphi _{\alpha
}^{\prime
}}\in \mathcal{J}(\langle\mathcal{M},e_\mathcal{N}\rangle)$ for every $%
\alpha\in \Lambda$ (see Lemma 1.2.1 of~\cite{Po2}). \\

Let $f_{\alpha }=\tau \circ \varphi _{\alpha }^{\prime }$. Then $\{f_{\alpha
}\}_{\alpha \in \Lambda }\subseteq \mathcal{M}_{\#}$. Since $\lim_{\alpha
}\Vert \varphi^{\prime }_{\alpha }(x)-x\Vert _{2,\tau }=0$ for every $x$ in $%
\mathcal{M}$, $\lim_{\alpha }f_{\alpha }(x)=\tau (x)$ for every
$x\in \mathcal{M}$.
Since $\mathcal{M} $ is the dual space of $\mathcal{M}_{\#}$, this implies that $%
\lim_{\alpha }f_{\alpha }=\tau $ in the weak topology on
$\mathcal{M}_\#$.
Since the weak closure and the strong closure of a convex set in $\mathcal{M}%
_{\#}$ are the same, $\tau $ is in the norm closure of the convex hull of $%
\{f_{\alpha }\}_{\alpha \in \Lambda }$. Note that $\tau \circ
(\sum_{i=1}^{n}\lambda _{\alpha _{i}}\varphi _{\alpha _{i}}^{\prime
})=\sum_{i=1}^{n}\lambda _{\alpha _{i}}f_{\alpha _{i}}.$ By taking
finitely many  convex combinations of $\{\varphi _{\alpha }^{\prime
}\}_{\alpha \in \Lambda }$, we can see that there exists a net
$\{\psi _{\beta}^{\prime
}\}_{\beta\in \Gamma}$ of completely positive $\mathcal{N}$-bimodular maps from $\mathcal{M}$
to $\mathcal{M}$ satisfying the conditions $2^{\prime }$ and $%
3^{\prime }$ in Theorem~\ref{Fang1}, $\psi _{\beta}^{\prime }(1)\leq
1$ for all $\beta\in \Gamma$ and the following condition:
$\lim_{\beta}\Vert g_{\beta}^{\prime }-\tau \Vert _{1}=0$ for
$g_{\beta}^{\prime}=\tau \circ \psi _{\beta}^{\prime }$.\newline

Let $b_\beta^{\prime }=1\vee (dg_\beta^{\prime }/d\tau )$ and $\psi
_{\beta}^{\prime\prime}(\cdot )=\psi _\beta ^{\prime }((b_\beta
^{\prime })^{-1/2}\cdot (b_\beta ^{\prime })^{-1/2})$. By
Lemma~\ref{Po2}, $\{\psi _\beta^{\prime\prime }\}_{\beta\in \Gamma}$
is a net of completely positive $\mathcal{N}$-bimodular maps from
$\mathcal{M}$
to $\mathcal{M}$, and satisfies $3^{\prime }$ in Theorem~%
\ref{Fang1}, $\psi _{\beta}^{\prime\prime }(1)\leq 1$ and $\tau\circ
\psi _{\beta}^{\prime\prime }\leq \tau$ for all $\beta\in \Gamma$.
By Lemma~\ref{Derivative}, $T_{\varphi _{\beta}^{\prime\prime
}}=T_{\psi _\beta^{\prime}}{b_\beta ^{\prime }}^{-1/2}J{b_\beta
^{\prime }}^{-1/2}J\in \langle\mathcal{M},e_\mathcal{N}\rangle$.
Since $T_{\psi _\beta^{\prime}}\in
\mathcal{J}(\langle\mathcal{M},e_\mathcal{N}\rangle)$,
 $T_{\psi _{\beta}^{\prime\prime }}\in \mathcal{%
J}(\langle\mathcal{M},e_\mathcal{N}\rangle)$ for every $\beta\in \Gamma$.
\newline

We may further assume that $\psi _{\beta}^{\prime\prime }(1)\leq
1-\epsilon_\beta$, $\epsilon_\beta>0$. Otherwise we can choose a net
of
positive numbers $\lambda_\beta$ with $0<\lambda_\beta<1$ and $%
\lim_\beta\lambda_\beta=1$ and consider $\lambda_\beta\cdot \psi
_{\beta}^{\prime\prime }$. Let $h_\beta=\psi _{\beta}^{\prime\prime }(1)$
and $k_\beta=d\tau\circ \psi _{\beta}^{\prime\prime }/d\tau$. By Lemma~\ref%
{equation}, there exist positive operators $a_\beta,b_\beta$ in $\mathcal{N}%
^{\prime }\cap \mathcal{M}$ such that $1-h_\beta=a_\beta E_{\mathcal{N}%
}(b_\beta)$ and $1-k_\beta=E_{\mathcal{N}}(a_\beta)b_\beta$.\newline

For every $\beta\in \Gamma$, define $\psi_\beta:\,
\mathcal{M}\rightarrow\mathcal{M}$ by
\begin{equation*}
\psi_\beta(x)=\psi _{\beta}^{\prime\prime}(x)+a_\beta E_{\mathcal{N}%
}(b_\beta x).
\end{equation*}
Clearly, every $\psi_\beta$ is a normal completely positive
$\mathcal{N}$-bimodular map. We have
\begin{equation*}
\psi_\beta(1)=\psi _{\beta}^{\prime\prime}(1)+a_\beta E_{\mathcal{N}%
}(b_\beta)=h_\beta+1-h_\beta=1,
\end{equation*}
and
\begin{equation*}
\tau(\psi_\beta(x))=\tau(\psi _{\beta}^{\prime\prime}(x))+\tau(a_\beta E_{%
\mathcal{N}}(b_\beta x))=\tau(xk_\beta)+\tau(E_{\mathcal{N}}(a_\beta)b_\beta
x)
\end{equation*}%
\begin{equation*}
=\tau(k_\beta x)+\tau((1-k_\beta)x)=\tau(x).
\end{equation*}
This proves that $\{\psi_\beta\}_{\beta\in \Gamma}$ satisfies the condition $%
1^{\prime }$ of Theorem~\ref{Fang1}.\newline

Note that $T_{\psi_\beta}=T_{\psi _{\beta}^{\prime\prime }}+a_\beta e_{%
\mathcal{N}} b_\beta$. Since $e_{\mathcal{N}}\in
\mathcal{J}(\langle\mathcal{M},e_\mathcal{N}\rangle)$, $T_{\psi
_{\beta}}\in \mathcal{J}(\langle\mathcal{M},e_\mathcal{N}\rangle)$
for every $\beta\in \Gamma$. This proves that
$\{\psi_\beta\}_{\beta\in \Gamma}$ satisfies the condition
$2^{\prime }$ of Theorem~\ref{Fang1}.\newline

Finally, for every positive operator $x$ in $\mathcal{M}$,
\begin{equation*}
\psi_\beta(x)-\psi _{\beta}^{\prime\prime}(x)=a_\beta E_{\mathcal{N}%
}(b_\beta x)\leq \|x\| a_\beta E_{\mathcal{N}}(b_\beta)=\|x\|(1-h_\beta)=\|x%
\|(1-\psi _{\beta}^{\prime\prime}(1)),
\end{equation*}
which shows that $\{\psi_\beta\}_{\beta\in \Gamma}$ satisfies the condition $%
3^{\prime }$ of Theorem~\ref{Fang1}.
\end{proof}

Let $\tau^{\prime }$ be another faithful normal trace on $\mathcal{M}$. Then
$\mathcal{M}\subseteq \mathcal{B}(L^2(\mathcal{M},\tau^{\prime }))$ is in
the standard form in the sense of Haagerup~\cite{Haa1}. Since the standard
representation of a von Neumann algebra is unique up to spatial isomorphism(see~%
\cite{Haa2}), the above arguments indeed prove the following
stronger result, which implies that the notion of ``relative
Haagerup property" considered by Boca in~\cite{Boca} is same as
Definition~\ref{HAP} given by Popa.

\begin{corollary}
Let $\mathcal{M}$ be a finite von Neumann algebra with a faithful normal trace $\tau$, and $\mathcal{N%
}$ a von Neumann subalgebra. Suppose $\{\varphi_\alpha\}_{\alpha\in \Lambda}$
is a net of normal completely positive $\mathcal{N}$-bimodular  maps from $\mathcal{M}$ to $%
\mathcal{M}$ satisfying the conditions $2$ and $3$ as in definition~\ref{HAP}%
, i.e. $\lim_{\alpha }\Vert \varphi _{\alpha }(x)-x\Vert _{2,\tau }=0$ for
all $x\in \mathcal{M}$ and the map $x\Omega\rightarrow
\varphi_\alpha(x)\Omega$ extends to a bounded operator $T_{\varphi_\alpha}$
in $\mathcal{J}(\langle\mathcal{M},e_\mathcal{N}\rangle)$ for every $%
\alpha\in \Lambda$. Then for every faithful normal trace $\tau'$ on $\mathcal{%
M}$, there exists a net $\{\psi_\beta\}_{\beta\in \Gamma}$ of normal completely positive $%
\mathcal{N}$-bimodular  maps from $\mathcal{M}$ to $\mathcal{M}$ satisfying

\begin{enumerate}
\item[$1^{\prime }$] $\psi_\beta(1)=1$ and $\tau'\circ \psi_\beta=\tau'$, $%
\forall \beta\in \Gamma$;

\item[$2^{\prime }$] $T_{\psi_\beta}\in \mathcal{J}(\langle\mathcal{M},e_%
\mathcal{N}\rangle)$,  $\forall \beta \in \Gamma$;

\item[$3^{\prime }$] $\lim_\beta\|\psi_\beta(x)-x\|_{2,\tau'}=0$, $\forall
x\in \mathcal{M}$.
\end{enumerate}

In particular, the relative Haagerup property does not depend on the
choice of faithful normal trace on $\mathcal{M}$.
\end{corollary}

\section{${\rm C_0}$-Correspondences}

We now show that Theorem \ref{Fang1} enables us to interpret
Haagerup's approximation property in the framework of Connes's
theory of correspondences. Suppose $\mathcal{M}$ is a finite von
Neumann algebra with a faithful normal trace $\tau$ and
$\mathcal{H}$ is a correspondence of $\mathcal{M}$.

\begin{definition}\label{C_0}
\emph{We say that $\mathcal{H}$ is a \emph{${\rm
C_0}$-correspondence} if $\mathcal{H}\cong\oplus_{\alpha\in
\Lambda}\mathcal{H}_{\varphi_{\alpha}}$, where each
$\mathcal{H}_{\varphi_{\alpha}}$ is the correspondence of
$\mathcal{M}$ associated to a completely positive map
$\varphi_\alpha:\, \mathcal{M}\rightarrow\mathcal{M}$  such that
 the extension operator $T_{\varphi_\alpha}$ of $\varphi_\alpha$
  is a compact operator in $\mathcal{B}(L^2(\mathcal{M},\tau))$.}
\end{definition}

\begin{remark}\emph{By the uniqueness of standard representation up to spatial
 isomorphism(see~\cite{Haa2}), the definition of
 ${\rm C_0}$-correspondence does not depend on the choice of
 $\tau$.}
\end{remark}

\begin{remark}\emph{The coarse correspondence $\mathcal{H}_{\rm co}$
of $\mathcal{M}$ is a ${\rm C_0}$-correspondence.
 By Proposition 1.2.5 of~\cite{Po2},
a sub-correspondence of a ${\rm C_0}$-correspondence (e.g., the
coarse correspondence) is not necessarily a ${\rm
C_0}$-correspondence. Let $\mathcal{H}_{\rm C_0}$ be the direct sum
of all  $\mathcal{H}_{\varphi}$ such that each
$\mathcal{H}_{\varphi}$ is the correspondence of $\mathcal{M}$
associated to a completely positive map $\varphi:\,
\mathcal{M}\rightarrow\mathcal{M}$ with
 the extension operator $T_{\varphi}$ of $\varphi$
  being a compact operator in $\mathcal{B}(L^2(\mathcal{M},\tau))$. Then
   $\mathcal{H}_{\rm C_0}$ is called the \emph{maximal ${\rm
   C_0}$-correspondence} of $\mathcal{M}$}.
\end{remark}

\begin{remark}\emph{Let $G$ be a discrete group, and  $\pi$ be a unitary
representation of $G$ on a Hilbert space $\mathcal{H}$. Then $\pi$
is unitarily equivalent to a direct sum of cyclic representations
$\pi_{f_\alpha}$ of $G$, where each $\pi_{f_\alpha}$ is the
representation associated to a positive definite function $f_\alpha$
on $G$. Recall that the representation $\pi$  is a ${\rm
C_0}$-representation if all matrix coefficients
$\omega_{\xi,\eta}(g)=\langle \pi(g)\xi,\eta\rangle$ belong to ${\rm
C_0}(G)$. It is easy to check that $\pi$ is a ${\rm
C_0}$-representation if and only if every $f_\alpha\in {\rm
C_0}(G)$. By~\cite{Haa1, Cho}, for every $f_\alpha$, there is a
unique normal completely positive map $\varphi_{f_\alpha}$ from the
group von Neumann algebra $L(G)$ to itself satisfying
$\varphi_{f_\alpha}(L_g)=f_\alpha(g)L_g$, where $L_g$ is the unitary
operator associated to $g$. By Lemma 1 and Lemma 2 of~\cite{Cho},
$f_\alpha$ is in ${\rm C_0}(G)$ if and only if the extension
operator $T_{\varphi_{f_\alpha}}$ of $\varphi_{f_\alpha}$ is a
compact operator in $\mathcal{B}(L^2(G))$. Hence, the correspondence
$\mathcal{H}_{\varphi_{f_\alpha}}$ of $L(G)$ associated to
$\varphi_{f_\alpha}$ is a ${\rm C_0}$-correspondence of
$\mathcal{M}$. So our definition of ${\rm C_0}$-correspondence of
finite von Neumann algebras is a natural analogue of the notion of
${\rm C_0}$-representation of groups.}
\end{remark}

The following theorem is the main result of this section.
\begin{theorem}
\label{HAPiffHaagerup}A finite von Neumann algebra
$(\mathcal{M},\tau)$ has  Haagerup's approximation property if and
only if the identity correspondence of $\mathcal{M}$ is weakly
contained in some
 ${\rm C_0}$-correspondence of $\mathcal{M}$.
\end{theorem}

To prove the above theorem, we need the following lemmas.\\

\begin{lemma}\label{Lemma:Coefficients} Let $\varphi$ be a normal completely positive map from
$\mathcal{M}$ to $\mathcal{M}$ such that
 the extension operator $T_{\varphi}$ of $\varphi$
  is a compact operator in $\mathcal{B}(L^2(\mathcal{M},\tau))$. Let $\xi=\sum_{i=1}^n a_i\otimes b_i$
be a vector in the  correspondence $\mathcal{H}_\varphi$ of
$\mathcal{M}$ associated to $\varphi$. Then $\xi$ is a left
$\tau$-bounded vector and the coefficient $\Phi_\xi$ corresponding
to $\xi$ is a normal completely positive map from $\mathcal{M}$ to
$\mathcal{M}$ such that $T_{\Phi_\xi}$ is a compact operator in
$\mathcal{B}(L^2(\mathcal{M},\tau))$.
\end{lemma}
\begin{proof} To see $\xi$ is a left $\tau$-bounded vector,
 we may assume that $\xi=a\otimes
b$. Then
\[\|\xi x\|_\varphi^2=\langle \xi x, \xi x\rangle_\varphi=\langle a\otimes (bx), a\otimes
(bx)\rangle_\varphi=\tau(\varphi(a^*a)bxx^*b^*)=\tau(x^*(b^*\varphi(a^*a)b)x)
\]
\[\leq \|b^*\varphi(a^*a)b\|\tau(x^*x)=
\|b^*\varphi(a^*a)b\|\|x\|_{2}^2.
\] Hence $\Phi_\xi$ is a normal
completely positive map from $\mathcal{M}$ to $\mathcal{M}$.
 For
every $x,y,z\in \mathcal{M}$, by equation (1) in section 2.6,
\[\langle\Phi_\xi(x)y\Omega, z\Omega\rangle_\tau=\langle x\xi y,\xi
z\rangle_\varphi=\sum_{i,j=1}^n\langle xa_j\otimes b_jy, a_i\otimes
b_iz\rangle_\varphi
\]
\[=\sum_{i,j=1}^n\tau(\varphi(a_i^*xa_j)b_jyz^*b_i^*)=\langle\sum_{i,j=1}^n b_i^*\varphi(a_i^*xa_j)b_j y\Omega,
z\Omega\rangle_\tau.
\]
This implies  that $\Phi_\xi(x)=\sum_{i,j=1}^n
b_i^*\varphi(a_i^*xa_j)b_j$. Hence, $\Phi_\xi$ can be extended to a
bounded operator from $L^2(\mathcal{M},\tau)$ to
$L^2(\mathcal{M},\tau)$ such that
\[T_{\Phi_\xi}=\sum_{i,j=1}^n b_i^*Jb_j^*JT_\varphi a_i^*Ja_j^*J.
\] Since $T_\varphi$ is a compact operator, $T_{\Phi_\xi}$ is also a compact
operator.
\end{proof}

\begin{lemma}\label{Lemma:convex cone}
Let $\mathcal{F}$ be the convex hull of the set of coefficients
$\Phi_\xi$ as in Lemma~\ref{Lemma:Coefficients}. Then $\mathcal{F}$
is a convex cone and for every $b\in \mathcal{M}$ and $\Phi\in
\mathcal{F}$, the completely positive map $b^*\Phi(\cdot) b$ belongs
to $\mathcal{F}$. Furthermore, $T_{\Phi}$ is a compact operator in
$\mathcal{B}(L^2(\mathcal{M},\tau))$ for all $\Phi\in \mathcal{F}$.
\end{lemma}
\begin{proof}  It is obvious that $\mathcal{F}$ is a convex cone.
To prove the rest, we may assume that $\Phi=\Phi_\xi$ is the
coefficient corresponding to $\xi\in \mathcal{H}_\varphi$ as in
Lemma~\ref{Lemma:Coefficients}.  Let $\eta=\xi
b=\sum_{i=1}^na_i\otimes b_ib\in \mathcal{H}$. By
Lemma~\ref{Lemma:Coefficients}, $\eta$ is a left $\tau$-bounded
vector. Let $\Phi_\eta$ be the coefficient corresponding to $\eta$.
By equation (1) in section 2.6,
\[\langle \Phi_\eta(x) y\Omega, z\Omega\rangle_{\tau}=\langle x\xi by, \xi bz\rangle_{\mathcal{\varphi}}
=\langle \Phi(x) by\Omega, bz\Omega\rangle_{\tau}=\langle
b^*\Phi(x)b y\Omega, z\Omega\rangle_{\tau}.\] This implies that
$\Phi_\eta=b^*\Phi b$. Hence $b^*\Phi b\in \mathcal{F}$.  By
Lemma~\ref{Lemma:Coefficients}, $T_{\Phi_\eta}$ is compact.
\end{proof}

\begin{lemma}\label{Lemma:sum of coefficients} Let $\mathcal{M}$ be a finite von Neumann
algebra with a faithful normal trace $\tau$, and $\mathcal{H}$,
$\mathcal{K}$ be two  correspondences of $\mathcal{M}$. Suppose
$\xi\in \mathcal{H}$ and $\eta\in \mathcal{K}$ are two left
$\tau$-bounded vectors, and $\Phi_{\xi}$, $\Phi_\eta$ are the
coefficients corresponding to $\xi$, $\eta$, respectively. Then
$\xi\oplus \eta$ is also a left $\tau$-bounded vector and
$\Phi_\xi+\Phi_\eta$ is the coefficient corresponding to $\xi\oplus
\eta\in \mathcal{H}\oplus\mathcal{K}$.
\end{lemma}
\begin{proof} It is clear that $\xi+\eta$ is a left $\tau$-bounded vector. By equation (1) in section 2.6,
\[\langle (\Phi_\xi+\Phi_\eta)(x)y\Omega,z\Omega\rangle_\tau=\langle x\xi y, \xi
z\rangle_{\mathcal{H}}+\langle x\eta y, \eta
z\rangle_{\mathcal{K}}=\langle x(\xi\oplus \eta) y, (\xi\oplus \eta)
z\rangle_{\mathcal{H}\oplus\mathcal{K}}\]\[=\langle
(\Phi_{\xi+\eta}(x)y\Omega,z\Omega\rangle_\tau.
\] Hence $\Phi_{\xi+\eta}=\Phi_\xi+\Phi_\eta$.
\end{proof}

Note that in the proof of Lemma 2.2 of~\cite{AdH}, if we replace the
arbitrary positive normal form $\phi$ (on line 10 of page 418) by an
arbitrary weak operator topology continuous positive form, then the
following lemma follows.

\begin{lemma}
\label{A-D Lemma}  Let $\Psi$ be a normal completely positive map
from $\mathcal{M}$ to $\mathcal{M}$. If $\Psi$ is in
 the closure of ${\mathcal F}$ in the pointwise weak operator topology. Then there exists a net $%
\{\Phi _{\alpha}\}_{\alpha\in\Lambda}$ in $\mathcal{F}$ such that $\Phi_{\alpha }(1)\leq \Psi (1)$ for all $%
\alpha\in\Lambda$, which converges to $\Psi$ in the pointwise weak
operator topology.
\end{lemma}

\begin{proof}[Proof of Theorem~\ref{HAPiffHaagerup}]
Suppose first that $\mathcal{M}$ has Haagerup's  approximation
property. By Theorem \ref{Fang1}, there is a net $(\varphi _{\alpha
})_{\alpha\in\Lambda}$ of unital normal completely positive maps
satisfying conditions ($1'$)-($3'$) in Theorem \ref{Fang1}.  It
immediately follows that the correspondence
$\mathcal{H}=\bigoplus_{\alpha\in\Lambda}
 \mathcal{H}_{\varphi _{\alpha }}$ is a ${\rm C_0}$-correspondence of
 $\mathcal{M}$ which weakly contains the identity correspondence of
 $\mathcal{M}$.\\

Conversely, suppose that $\mathcal{H}$ is a ${\rm
C_0}$-correspondence of $\mathcal{M}$ which weakly contains the
identity correspondence of
 $\mathcal{M}$.
We may assume $\mathcal{H}=\oplus_{\beta\in
\Gamma}\mathcal{H}_{\varphi_\beta}$, with each $\varphi_\beta:\,
\mathcal{M}\rightarrow\mathcal{M}$ is a normal completely positive
map  such that the extension operator $T_{\varphi_\beta}$ of
$\varphi_\beta$
 is a compact operator in $\mathcal{B}(L^2(\mathcal{M},\tau))$.
 Since the identity correspondence of $\mathcal{M}$ is weakly contained in
  $\mathcal{H}$, for every
$\epsilon>0$ and every finite subset $E$ of  $\mathcal{M}$, there
exists a $\xi\in \mathcal{H}$ such that
\[|\langle x\xi y,\xi z\rangle_{\mathcal{H}}-\langle x\Omega y,\Omega z\rangle_\tau|<\epsilon,\quad \forall x,y,z\in
E.
\]
We may assume that $\xi=\xi_1\oplus \cdots\oplus \xi_n$, where
$\xi_i=\sum_{j=1}^{n_i}a_{ij}\otimes b_{ij}\in
\mathcal{H}_{\varphi_{\beta_i}}$. Let $\Phi_\xi$ be the coefficient
corresponding to $\xi$. By Lemma~\ref{Lemma:convex cone} and
Lemma~\ref{Lemma:sum of coefficients}, $\Phi_\xi\in \mathcal{F}$. By
equation (1) in section 2.6,
\[|\langle \Phi_\xi(x)y\Omega,z\Omega\rangle_\tau-\langle x\Omega y,\Omega
z\rangle_\tau|=|\langle x\xi y,\xi z\rangle_{\mathcal{H}}-\langle
x\Omega y,\Omega z\rangle_\tau|<\epsilon, \quad \forall x,y,z\in E.
\] This implies that  there exists a net
$(\Phi_{\alpha' })_{\alpha'\in \Lambda'}$ of completely positive
maps in $\mathcal{F}$ such that $\lim_{\alpha' }\Phi_{\alpha'
}(x)=x$ in the weak operator topology for every $x\in\mathcal{M}$.\\

 By Lemma~\ref{A-D Lemma}, there is a net $\{\varphi_\alpha\}_{\alpha\in \Lambda}$ in $\mathcal{F}$ such
that $\lim_{\alpha}\varphi_{\alpha }(x)=x$ in the weak operator
topology for every $x\in \mathcal{M}$ and   $\varphi _{\alpha
}(1)\leq 1$ for every $\alpha\in \Lambda$.
 Now given $x\in \mathcal{M}$:%
\begin{eqnarray*}
||\varphi _{\alpha }(x)-x||_{2} &=&\tau (\varphi _{\alpha }(x)^{\ast
}\varphi _{\alpha }(x))+\tau (x^{\ast }x)-2\mathrm{Re}\, \tau (x^{\ast
}\varphi _{\alpha }(x)) \\
&\leq &||\varphi _{\alpha }(1)||\tau (\varphi _{\alpha }(x^{\ast }x))+\tau
(x^{\ast }x)-2\mathrm{Re}\,\tau (x^{\ast }\varphi _{\alpha }(x)) \\
&\leq &\tau (\varphi _{\alpha }(x^{\ast }x))+\tau (x^{\ast }x)-2\mathrm{Re}%
\,\tau (x^{\ast }\varphi _{\alpha }(x))\text{.}
\end{eqnarray*}%
Since $\lim_{\alpha }\varphi _{\alpha }(x)=x$ in the weak operator
topology for every $x\in \mathcal{M}$ it follows that $\lim_{\alpha
}\tau (\varphi _{\alpha }(x^{\ast }x))=\tau (x^{\ast }x)$ and
$\lim_{\alpha }\tau (x^{\ast }\varphi _{\alpha }(x))=\tau (x^{\ast
}x)$. Therefore $\lim_{\alpha }||\varphi _{\alpha }(x)-x||_{2}=0$.
This proves that $(\varphi _{\alpha })_{\alpha }$ is a net of
completely positive maps that approximate the identity pointwise in
the trace-norm. Since $\varphi_\alpha\in \mathcal{F}$,
it follows that $T_{\varphi_\alpha}$ is a compact operator on $L^2(\mathcal{M},\tau)$. By Theorem \ref%
{Fang1}, $\mathcal{M}$ has Haagerup's approximation property.
\end{proof}

As an application of Theorem~\ref{HAPiffHaagerup}, we prove the
following theorem.

\begin{theorem}If
  $\mathcal{M}$ has Haagerup's approximation property, then the class of ${\rm
 C_0}$-correspondences of $\mathcal{M}$ is dense in ${\rm
 Corr}(\mathcal{M})$.
\end{theorem}
\begin{proof} By section 2.6, it is clear that we need only to prove
that every cyclic correspondence $\mathcal{H}_\Phi$ of $\mathcal{M}$
associated to a coefficient $\Phi$ belongs to the closure of the set
of ${\rm C_0}$-correspondences of $\mathcal{M}$. Since $\mathcal{M}$
has Haagerup's approximation property, there is a net
$\{\varphi_\alpha\}_{\alpha\in\Lambda}$ of normal completely
positive maps of $\mathcal{M}$, such that
\begin{enumerate}
\item  $\varphi_\alpha(1)=1$, $\forall \alpha\in \Lambda$,
 \item  $T_{\varphi_\alpha}$ is compact, $\forall \alpha\in \Lambda$,
\item $\lim_\alpha\|\varphi_\alpha(x)-x\|_2=0$, $\forall x\in
\mathcal{M}$.
\end{enumerate}
 Hence,  each $T_{\Phi\circ\varphi_\alpha}=T_\Phi
T_{\varphi_\alpha}$ is compact and
$\lim_\alpha\|\Phi\circ\varphi_\alpha(x)-\Phi(x)\|_2=0$ for every
$x\in \mathcal{M}$. By Remark 2.1.4 of~\cite{Po1},
$\mathcal{H}_{\Phi\circ\varphi_\alpha}\rightarrow \mathcal{H}_\Phi$.
\end{proof}

\begin{corollary}If
  $\mathcal{M}$ has Haagerup's approximation property, then every correspondence of $\mathcal{M}$
  is weakly contained in $\mathcal{H}_{\rm C_0}$, the maximal ${\rm C_0}$-correspondence of $\mathcal{M}$.
\end{corollary}

\section{Relative Amenability and Haagerup's Approximation Property}

Popa asks the following question in~\cite{Po2} (see Section 3.5.2): If $\mathcal{%
N\subseteq M}$ is an inclusion of finite von Neumann algebras and $\mathcal{N%
}$ has Haagerup's approximation property, does $\mathcal{M}$ also
have Haagerup's approximation property? The following theorem
answers Popa's question affirmatively.

\begin{theorem}
\label{RelAmen}If $\mathcal{N}\subseteq \mathcal{M}$ is an amenable
inclusion of finite von Neumann algebras and $\mathcal{N}$ has
Haagerup's approximation property then $\mathcal{M}$ also has
Haagerup's approximation property.
\end{theorem}

To prove Theorem~\ref{RelAmen}, we need the following lemmas.

\begin{lemma}
\label{CyclicInduction}Let $\mathcal{N\subseteq M}$ be an inclusion
of finite von Neumann algebras, and $E_{\mathcal{N}}$ be the normal
$\tau $-preserving conditional expectation of $\mathcal{M}$ onto
$\mathcal{N}$. If $\mathcal{H}_{\varphi }$ is the  correspondence of
$\mathcal{N}$ associated to a normal completely positive map
$\varphi$ from $\mathcal{N}$ to $\mathcal{N}$ and
$\mathcal{H}_{\varphi \circ E_{\mathcal{N}}}$ is the correspondence
of $\mathcal{M}$  associated to the normal completely positive map
$\varphi\circ E_{\mathcal{N}}$ from $\mathcal{M}$ to $\mathcal{M}$,
then ${\rm
Ind}_{\mathcal{N}}^{\mathcal{M}}(\mathcal{H}_\varphi)\cong
\mathcal{H}_{\varphi \circ E_{\mathcal{N}}}$.
\end{lemma}

\begin{proof} Denote by
$\mathcal{K}={\rm
Ind}_{\mathcal{N}}^{\mathcal{M}}(\mathcal{H})=L^2(\mathcal{M})
\underset{\mathcal{N}}{\otimes}\mathcal{H}_\varphi\underset{\mathcal{N}}{\otimes}
L^2(\mathcal{M})$, where the first $L^2(\mathcal{M})$ is regarded as
a left $\mathcal{M}$ and right $\mathcal{N}$ module and the second
$L^2(\mathcal{M})$ is regarded as a left $\mathcal{N}$ and right
$\mathcal{M}$ module. Let $\xi\in \mathcal{H}_\varphi$ be the vector
corresponding to $\Omega\otimes \Omega$, which is a cyclic vector of
$\mathcal{H}_\varphi$. Given $x_{1}, x_{2}, y_1, y_2\in
\mathcal{M}$, we have
\[\langle (x_{1}\otimes (\xi \otimes y_{1}),\,x_{2}\otimes (\xi
\otimes y_{2})\rangle _{\mathcal{K}}  =\langle q(\xi \otimes
y_{1}),\, \xi \otimes y_{2}\rangle
_{\mathcal{H}_\varphi\underset{\mathcal{N}}{\otimes}
L^2(\mathcal{M})},
\] where
$q\in\mathcal{N}$ is the Radon-Nikodym derivative of $\mathcal{N}\ni
z\mapsto \langle x_1z,\, x_2 \rangle_{L^2(\mathcal{M})}$ with
respect to $\tau_{\mathcal{N}}$. Note that
\[\langle x_1z,\, x_2
\rangle_{L^2(\mathcal{M})}=\tau(zx_2^*x_1)=\tau(zE_{\mathcal{N}}(x_2^*x_1)).
\]
Hence $q=E_{\mathcal{N}}(x_2^*x_1)$ and
\[\langle (x_{1}\otimes (\xi \otimes y_{1}),\,x_{2}\otimes (\xi
\otimes y_{2})\rangle _{\mathcal{K}} =\langle
E_{\mathcal{N}}(x_2^*x_1)\xi \otimes y_{1},\, \xi \otimes
y_{2}\rangle _{\mathcal{H}_\varphi\underset{\mathcal{N}}{\otimes}
L^2(\mathcal{M})}=\langle E_{\mathcal{N}}(x_2^*x_1)\xi p,\, \xi
\rangle _{\mathcal{H}_\varphi},
\]
where $p\in\mathcal{N}$ is the Radon-Nikodym derivative of
$\mathcal{N}\ni z\mapsto \langle zy_1,\, y_2
\rangle_{L^2(\mathcal{M})}$ with respect to $\tau_{\mathcal{N}}$.
Note that
\[\langle zy_1,\, y_2
\rangle_{L^2(\mathcal{M})}=\tau(zy_1y_2^*)=\tau(zE_{\mathcal{N}}(y_1y_2^*)).
\]
Hence $p=E_{\mathcal{N}}(y_1y_2^*)$ and
\[\langle (x_{1}\otimes (\xi \otimes y_{1}),\,x_{2}\otimes (\xi
\otimes y_{2})\rangle _{\mathcal{K}} =\langle
E_{\mathcal{N}}(x_2^*x_1)\xi p,\, \xi \rangle
_{\mathcal{H}_\varphi}=\langle E_{\mathcal{N}}(x_2^*x_1)\xi
E_{\mathcal{N}}(y_1y_2^*),\, \xi \rangle _{\mathcal{H}_\varphi}
\]
\[=\tau(\varphi(E_{\mathcal{N}}(x_2^*x_1))E_{\mathcal{N}}(y_1y_2^*))=\tau(\varphi(
E_{\mathcal{N}}(x_2^*x_1))y_1y_2^*)=\langle x_{1}\xi  y_{1},\,x_{2}
\xi y_{2}\rangle _{\mathcal{H}_{\varphi\circ E_{\mathcal{N}}}}.
\]

 Therefore the map defined on simple
tensors by $(x_{1}\otimes \xi )\otimes
x_{2}\mapsto x_{1}\xi x_{2}$ extends to an $\mathcal{M}$-linear isometry from $%
{\rm Ind}_{\mathcal{N}}^{\mathcal{M}}(\mathcal{H}_\varphi)$ onto
$\mathcal{H}_{\varphi \circ E_{\mathcal{N}}}$.
\end{proof}

The proof of the following lemma is an easy exercise.
\begin{lemma}
\label{DirectSumInduction} Let $\mathcal{N\subseteq M}$ be an
inclusion of finite von Neumann algebras, and $E_{\mathcal{N}}$ be
the normal $\tau $-preserving conditional expectation of
$\mathcal{M}$ onto $\mathcal{N}$. Suppose for $\alpha\in \Lambda$,
$\mathcal{H}_{\varphi_\alpha}$ is the  correspondence of
$\mathcal{N}$ associated to a normal completely positive map
$\varphi_\alpha$ from $\mathcal{N}$ to $\mathcal{N}$.  Then ${\rm
Ind}_{\mathcal{N}}^{\mathcal{M}}(\oplus_{\alpha\in
\Lambda}\mathcal{H}_{\varphi_\alpha})\cong \oplus_{\alpha\in
\Lambda}\mathcal{H}_{\varphi_\alpha \circ E_{\mathcal{N}}}$.
\end{lemma}

\begin{lemma}
\label{C0Lemma}If $\mathcal{N\subseteq M}$ is an inclusion of finite
von Neumann algebras and $\mathcal{H}$ is a ${\rm
C_{0}}$-correspondence of $\mathcal{N}$, then ${\rm
Ind}_{\mathcal{N}}^{\mathcal{M}}(\mathcal{H})$ is a ${\rm
C_{0}}$-correspondence of $\mathcal{M}$.
\end{lemma}

\begin{proof} Let $E_{\mathcal{N}}$ be the
$\tau $-preserving normal conditional expectation of $\mathcal{M}$ onto $%
\mathcal{N}$. Suppose $\mathcal{H}=\bigoplus\nolimits_{\alpha \in
\Lambda}\mathcal{H}_{\varphi _{\alpha }} $  such that
 $T_{\varphi_\alpha}$ is a compact operator in
$\mathcal{B}(L^2(\mathcal{N},\tau))$.
 By Lemma \ref{CyclicInduction} and Lemma \ref%
{DirectSumInduction} we have that ${\rm Ind}_{\mathcal{N}}^{\mathcal{M}%
}(\bigoplus\nolimits_{\alpha \in \Lambda}\mathcal{H}_{\varphi
_{\alpha }})\cong \bigoplus\nolimits_{\alpha \in
\Lambda}\mathcal{H}_{\varphi _{\alpha }\circ E_{\mathcal{N}}}$.
 Since $T_{\varphi _{\alpha }\circ E_{\mathcal{N}}}=T_{\varphi _{\alpha }}
e_{\mathcal{N}}$, the operator $T_{\varphi _{\alpha }\circ
E_{\mathcal{N}}}$ is a compact operator in
$\mathcal{B}(L^2(\mathcal{M},\tau))$.
 So ${\rm
Ind}_{\mathcal{N}}^{\mathcal{M}}(\mathcal{H})$ is a ${\rm
C_{0}}$-correspondence of $\mathcal{M}$.
\end{proof}

\begin{proof}[Proof of Theorem~\ref{RelAmen}]
Let $\mathcal{H}$ be a ${\rm C_0}$-correspondence of $\mathcal{N}$
that weakly contains the identity correspondence
$L^2(\mathcal{N},\tau)$ of
 $\mathcal{N}$. By
Lemma~\ref{C0Lemma}, ${\rm
Ind}_{\mathcal{N}}^{\mathcal{M}}(\mathcal{H})$ is a ${\rm
C_{0}}$-correspondence of $\mathcal{M}$. Note that
$L^2(\mathcal{N},\tau)\prec \mathcal{H}$. By the continuity of
induction operation (see Proposition 2.2.1 of~\cite{Po1}),  we see
that
\begin{equation*}
{\rm Ind}_{\mathcal{N}}^{\mathcal{M}}(L^2(\mathcal{N},\tau))\prec
{\rm Ind}_{\mathcal{N}}^{\mathcal{M}}(\mathcal{H})\text{.}
\end{equation*}%
Since $\mathcal{N}\subseteq \mathcal{M}$ is an amenable inclusion, we have%
\begin{equation*}
L^2(\mathcal{M},\tau)\prec \mathcal{H}_{\mathcal{N}}={\rm Ind}_{\mathcal{N}}^{\mathcal{M}%
}(L^2(\mathcal{N},\tau)).
\end{equation*}%
By the transitivity of $\prec $ we obtain%
\begin{equation*}
L^2(\mathcal{M},\tau)\prec {\rm
Ind}_{\mathcal{N}}^{\mathcal{M}}(\mathcal{H}) \text{.}
\end{equation*}
 By Theorem~\ref{HAPiffHaagerup}, $\mathcal{M}$ has
Haagerup's approximation property.
\end{proof}

\end{document}